\documentclass[11pt]{article}
\usepackage{amssymb}
\usepackage{latexsym}
\begin{document}
\setlength{\baselineskip}{15pt}
\title{Transition between extinction and blow-up \\in a generalized Fisher-KPP model}
\author{ Benito Hern\'{a}ndez-Bermejo $^{a,1}$ \mbox{} \mbox{} \mbox{} Ariel S\'{a}nchez-Vald\'{e}s $^b$
}
\date{}

\maketitle
\begin{center}
{\em 
$^a$ Departamento de F\'{\i}sica. Universidad Rey Juan Carlos. Calle Tulip\'{a}n S/N. 
28933--M\'{o}stoles--Madrid. Spain. \\ \mbox{} \\
$^b$ Departamento de Matem\'{a}tica Aplicada. Universidad Rey Juan Carlos. 
Calle Tulip\'{a}n S/N. 28933--M\'{o}stoles--Madrid. Spain. } 
\end{center}

\mbox{}

\mbox{}

\mbox{}

\mbox{}

\begin{center} 
{\bf Abstract}
\end{center}
\noindent
Stationary solutions of the Fisher-KPP equation with general nonlinear diffusion and arbitrary reactional kinetic orders terms are characterized. Such stationary (separatrix-like) solutions disjoint the blow-up solutions from those showing extinction. In addition a criterion for general parameter values is presented, which allows determining the blow-up or vanishing character of the solutions.

\mbox{}

\mbox{}

\mbox{}


\mbox{}

\mbox{}



\noindent {\bf PACS (2010) codes:} 02.30.Jr, 47.20Ky, 82.40.Bj.

\vfill

\noindent $^1$ Corresponding author. Telephone: (+34) 91 488 73 91. Fax: (+34) 91 664 74 55. \newline 
\mbox{} \hspace{0.05cm} E-mail: benito.hernandez@urjc.es 

\pagebreak
\begin{flushleft}
{\bf 1. Introduction}
\end{flushleft}

\

The following reaction-diffusion equation is to be considered in what follows:
\begin{equation}
u_t= \kappa (u^{m-1}u_x)_{x}+ \alpha u^p- \beta u^q \; , \qquad x \in \mathbb{R} \; , \; \; t>0 \; ,
\label{ec0}
\end{equation}
with $\alpha , \beta , \kappa , m>0$ and $p,q \in \mathbb{R}$. After the change of variables $x \to ax$, $t\to bt$, $u\to lu$, where $a=( \kappa l^{m-p}/\alpha)^{1/2}$, $b=(l^{1-p}/\alpha)$, and $l=(\beta/\alpha)^{1/(p-q)}$, equation (\ref{ec0}) becomes:
\begin{equation}
u_t=(u^{m-1}u_x)_{x}+ u^p- u^q \; , \qquad x \in \mathbb{R} \; , \; \; t>0 \; .
\label{ec1}
\end{equation}

Equation (\ref{ec1}) has been analyzed in detail in the literature. For
instance, equation (\ref{ec1}) with $m=1$, $p=1$, and $q=2$ was proposed by Fisher \cite{fis}
and Kolmogorov {\em et al.} \cite{kpp}. In addition, for the more general situation $p\leq q$, $m>0$ and $m+p>0$ the reader is referred to \cite{arw}-\cite{SH}. The complementary case $p>q$ will be of interest in this work and has been previously considered by several authors, specially in the search of travelling waves and related applications \cite{GK},\cite{SH}-\cite{egs}. It is worth emphasizing in this context that the consideration of general reaction exponents $p$ and $q$ such as those used in (\ref{ec1}) is of central interest in the framework of applied modelling (for instance, see \cite{VO,M} and references therein for some classical applications). The same equation (either for $p \leq q$ or $p > q$) is also relevant from the applied perspective for the modelling of growth and diffusion processes (see \cite{banks} for a classical presentation, as well as \cite{vos}-\cite{marus2} for specific developments including the case $p>q$).

Another two particular cases considered for equation (\ref{ec0}) correspond to either $\beta=0$ or $\alpha=0$. In such situations we have, respectively, a reaction-diffusion equation
\begin{equation}
\label{r1}
	u_t=(u^{m-1}u_x)_{x}+\alpha u^p, \qquad x \in \mathbb{R} \; , \; \; t>0 \; ,
\end{equation}
or an absorption-diffusion equation:
\begin{equation}
\label{r2}
	u_t=(u^{m-1}u_x)_{x}-\beta u^q, \qquad x \in \mathbb{R} \; , \; \; t>0  \; .
\end{equation}
There exist a numerous literature dealing with the previous two equations (\ref{r1}) and 
(\ref{r2}) and their physical implications, for example see \cite{ker}-\cite{GV3}. Some of the most significant features of these equations are behaviors such as blow-up for (\ref{r1}) if $p>1$ (in which the solution tends to infinity in a finite time, for instance see \cite{sgk,DL,GV3,VH}) and the extinction for (\ref{r2}) provided $q<1$ (in which the solution tends to zero in finite time, e.g. see \cite{K2}-\cite{FV}). 

The main result of the present work is that the previous solutions do not exhaust all the possibilities for that situation. To be precise, in what follows we shall consider the following Cauchy problem:
\begin{equation}
\label{p1}
\left\{
\begin{array}{ll} 
u_t=(u^{m-1}u_x)_{x}+ u^p- u^q \: , & x\in \mathbb{R} \: , \: t>0 \: , \: 
	p > q > 0 \: , \: m > 1  \; ,  \\
u(x,0)=u_0(x) \: , & x\in\mathbb{R} \; , u_0(x) \geq 0 \; .
\end{array}
\right. 
\end{equation}

The initial condition $u_0(x) \in L_1 \cap C$ (integrable and continuous functions). For instance, if $u(x)$ represents a chemical concentration  (number of molecules per unit lenght) then the condition $\int_{-\infty}^{\infty} u_0(x) \: \mbox{d}x < \infty$ accounts for the (finite) initial number of molecules. Regarding the reactional structure associated with the previous system (as described by the purely reactional o.d.e. $u_t = \alpha u^p-\beta u^q$) it is worth recalling that it can be easily determined from simple reactional schemes according to the standard stoichiometric procedures \cite{eth} for both cases $p>q$ and $p<q$. Recall also that there is no restriction in the values of $p$ and $q$ for the growth with diffusion phenomena, as displayed in \cite{vos}-\cite{marus2}.

Before going into the main contributions of the work, some previous results and terminology are required: 

\mbox{}

\noindent{\bf Definition 1.} 
{\em
\begin{enumerate}
  \item We call $\underline{u}$ a subsolution of (\ref{ec1}) if it satisfies
\begin{equation}
 \label{p3}
 \left\{
 \begin{array}{ll}
 \underline{u}_t\leq (\underline{u}^{m-1}\underline{u}_x)_{x}+\underline{u}^p-\underline{u}^q
 \: , & x\in \mathbb{R} \: , \: t>0  \; , \\ 
 \underline{u}(x,0)\leq u_0(x) \: , & x\in\mathbb{R}  \; .
 \end{array}
 \right. 
\end{equation}
  \item We call $\overline{u}$ a supersolution of (\ref{ec1}) if it satisfies (\ref{p3}) with opposite
  inequalities.
\end{enumerate}
}

\mbox{}

The following result, which is classical for this type of equations \cite{ker,sgk,bor,Ka1}, shall be employed along the present work:

\mbox{}

\noindent{\bf Proposition 1.} 
{\em 
Let  $u_1(x,t)$, $u_2(x,t)$ be two sub (correspondingly, super) solutions for every $0<t<T$, then:
$$
u_1(x,0)\le u_2(x,0) \:\:\; \mbox{in } \mathbb{R} \;\Rightarrow\; 
u_1(x,t)\le u_2(x,t) \:\:\; \mbox{in } \mathbb{R}\times(0,T)  \; .
$$
}

\mbox{}

In order to fix the terminology, the following definition is also necessary: 

\mbox{}

\noindent{\bf Definition 2.} 
{\em Let $T$ be the maximal time of existence of solution $u(x,t)$. If $T<\infty$ and  
$\lim\limits_{t\to T}\|u(x,t)\|_\infty=\infty$, then this is termed blow-up. If
$T<\infty$ and $\lim\limits_{t\to T}\|u(x,t)\|_\infty=0$ then this is termed extinction. Finally, 
if $T=\infty$ we have growth if $\lim\limits_{t\to \infty}\|u(x,t)\|_\infty = \infty$. 
}

\mbox{}

The main results can now be developed. The structure of this work is the following. In Section 2, the stationary solutions of (\ref{p1}) are characterized. Section 3 is devoted to show that for $p \geq m >q$, $p \geq 1$, $q \leq 1$, there exists an stationary solution (separatrix-like solution) which disjoints the increasing solutions (presenting blow-up if $p>1$, and growth if $p=1$) from the decreasing ones (showing extinction if $q<1$, and vanishing at infinity if $q=1$). In Section 4, for the more general case $p > q > 0$, $p>1$, $m>0$, we shall prove that there exist initial conditions for which the solution blows-up as well as initial conditions for which the solution presents extinction. The work is concluded in Section 5 with some examples, including the construction of a new explicit solution. 

\mbox{}

\mbox{}

\pagebreak
\begin{flushleft}
{\bf 2. Stationary solutions}
\end{flushleft}

The following result constitutes the basis of the results to be developed in the next section: 

\mbox{}

\noindent{\bf Proposition 2.}
{\em The stationary solutions of equation (\ref{ec1}) in the case $m+q>0$ are given by:
\begin{equation}
\label{cons}
\pm k_1 x+C=2 _2F_1(1/2,k_2;3/2;1-g(x))\sqrt{1-g(x)}  \; ,
\end{equation}
where $_2F_1$ is the hypergeometric function, $k_1=(p-q)\sqrt{\frac{2}{m+q}}(\frac{m+p}{m+q})^{\frac{q-m}{2(p-q)}}$, and $k_2=1+\frac{q-m}{2(p-q)}$.
}

\mbox{}

\noindent{\bf Proof.}

Let us recall that according to S\'{a}nchez-Vald\'{e}s and Hern\'{a}ndez-Bermejo \cite{SH}, for $m+q>0$ there exists a stationary solution $u(x,t)=f(x)$ of problem (\ref{ec1}), and such solution verifies the differential equation:
\begin{equation}
\label{edbase}
	(f^{m-2}f')^2-\frac{2f^{m+q-2}}{m+q}+\frac{2f^{m+p-2}}{m+p}=0  \; .
\end{equation}
Now making use of the transformation $g = \frac{m+q}{m+p}f^{p-q}$ we arrive at:
\begin{equation}
\label{g}
    g'=\pm k_1g^{k_2}\sqrt{1-g}  \; ,
\end{equation}
where $k_1$ and $k_2$ have the expressions given in the statement of the proposition. If we now integrate equation (\ref{g}) we obtain the solution as the implicit closed-form (\ref{cons}).
$\:\;\:\; \Box$

\mbox{}

\textbf{Remark.} From equation (\ref{g}) it can be seen that:
\begin{itemize}
\item  $g=1$ is a maximum of $g(x)$, and it corresponds to $f=(\frac{m+p}{m+q})^{1/(p-q)}>1$.
\item Note that if $g_s(x)$ is a solution of (\ref{g}), then all the solutions of (\ref{g}) are translations of $g_s$ in $x$, namely they are of the form $g_s(x+k)$ for every $k \in \mathbb{R}$.  
\item If $m \leq q$ then solution $g(x)$, and therefore also $f(x)$, are strictly positive.
\item If $m > q$ then solution $g(x)$, and therefore also $f(x)$, have compact support, namely they are strictly positive in a bounded interval, and they vanish in the complementary of such bounded interval. 
\item  $k_2=1$ ($m=q$) is a limit case separating the compact support solutions from those 	that tend to zero in infinity. In addition, in the case in which compact support exists, 	$k_2<1$, the size of the support can be computed as: 
\[
2\int_0^1 \frac{dg}{k_1g^{k_2}\sqrt{1-g}}= \frac{2}{k_1} B \left( 1-k_2,\frac12 \right)  \; ,
\] 
where $B(a,b)$ is the Beta function.	 
\item For the properties of the hypergeometric and beta functions, the reader is referred to \cite{abr,arf}.
\item For some values of $k_2$ explicit-form solutions can be found. In the following examples the solution is centered without loss of generality at the point $x=0$:
\begin{description}
  \item[{\em (a).}] $k_2=0$ ($p=\frac{m+q}{2}$): $g(x)=1-(k_1x/2)^2$.
  \item[{\em (b).}] $k_2=1/2$ ($p=m$): $g(x)=(\cos(k_1x)+1)/2$.
  \item[{\em (c).}] $k_2=1$ ($q=m$): $g(x)=1-(\tanh(k_1x/2))^2$.
  \item[{\em (d).}] $k_2=3/2$ ($q=\frac{m+p}{2}$): $g(x)=4/((k_1 x)^2+4)$.
\end{description}
\end{itemize}

\mbox{}

\mbox{}

\begin{flushleft}
{\bf 3. The stationary solution as a separatrix}
\end{flushleft}

In what follows, the stationary solutions considered in Section 2 are characterized and employed in order to develop a criterion for the behavior of the solutions of (\ref{p1}): 

\mbox{}

\noindent{\bf Theorem 1.}
{\em Let $u_0(x) \in L_1 \cap C$, $u_0(x) \geq 0$, be an intitial condition of (\ref{ec1}) such that $p \geq m > q$, $p \geq 1$, $q \leq 1$. And let $E(x)$ be a stationary solution verifying (\ref{edbase}). We then have: 
\begin{enumerate}
    \item If $u_0(x)>E(x)$ then $u(x,t)$ blows-up in finite time for $p>1$ and $u(x,t)$ presents growth if $p=1$. 
    \item If $u_0(x)<E(x)$ then $u(x,t)$ presents extinction in finite time for $q<1$, and $u(x,t)$ vanishes in infinite 			time for $q=1$.
\end{enumerate}
}

\mbox{}

\noindent{\bf Proof.}

In order to prove the first statement of the Theorem, let us first show that $G(x,t)=(T-at)^{-\alpha} E(x)$, is a subsolution of (\ref{ec1}) provided $p>1$ and $T<1$. Therefore, the following inequality must be proved: 
\begin{equation}
\label{llator}
a\alpha(T-at)^{-\alpha-1}E\leq (T-at)^{-\alpha m}(E^{m-1}E')'+(T-at)^{-\alpha p}E^p -(T-at)^{-\alpha q}E^q  \; .
\end{equation}
In the points in which $E=0$, equation (\ref{llator}) is identically verified. Thus the proof is continued for the points in which it is $E>0$ (the support of $E$). 
Taking into account that $E$ is a stationary solution of (\ref{ec1}) as described in Section 2, namely $(E^{m-1}E')'=E^q-E^p$, we have that:
\[
a\alpha(T-at)^{-\alpha-1}E\leq E^q ((T-at)^{-\alpha m}-(T-at)^{-\alpha q})+
E^p ((T-at)^{-\alpha p}-(T-at)^{-\alpha m})  \; .
\]
Defining $k= \min\{1-T^{\alpha(m-q)},1-T^{\alpha(p-m)}\}$ for the case $p>m$, since $(T-at)<T$ it suffices to prove that:
\[
a\alpha/k\leq E^{q-1} (T-at)^{-\alpha m+\alpha+1}+E^{p-1}(T-at)^{-\alpha p+\alpha+1}=F(E)  \; .
\]
Since $F(E)>0$ and $F(E)\to \infty$ when $E\to 0$ and $E\to \infty$, then $F(E)$ reaches a minimum for every fixed $t$. Let us show that such minimum value has a lower bound by a $C>0$ for every $t$. For $q<1$ we have that
\[
F'(E)=(q-1)E^{q-2}(T-at)^{-\alpha m+\alpha+1}+(p-1)E^{p-2}(T-at)^{-\alpha p+\alpha+1}  \; ,
\] 
and thus $F'(E)=0$ if
\[
E= \left( \frac{1-q}{p-1}(T-at)^{\alpha(p-m)} \right)^{1/(p-q)}=E_0  \; .
\]
Setting $\alpha=(p-q)(m-q)^{-1}(p-1)^{-1}$ we have that $F(E_0)=\frac{p-q}{p-1}(\frac{1-q}{p-1})^{(p-1)/(p-q)}$. In addition, if $q=1$ we choose $\alpha=(m-1)^{-1}$, and the minimum of $F(E)$ is achieved when $E=0$ and $F(0)=1$. 
Consequently, taking a sufficiently small $a$ inequality (\ref{llator}) is proved. For the case $p=m$ we should take $k=1-T^{\alpha (m-q)}$, with $\alpha = 1/(m-1)$ and the proof is direct after taking a sufficiently small $a$.

Now if we set $T=(\max_{x\in S}\{E(x)/u_0(x)\})^{1/\alpha}$, where $S=\{x:E(x) > 0\}$, then we have that $E(x)\leq T^{-\alpha}E(x)\leq u_0(x)$, and according to Proposition 1 this implies $G(x)\leq u(x,t)$. Taking into account that $G(x)$ blows-up, then so does $u(x,t)$.

Additionally, if $p=1$ then we set $G(x,t)=(T+ at)^{\alpha}E(x)$ with $T>1$, and proceeding in a way analogous to that of the $p>1$ case, the outcome is: 
\[
		a \alpha /k \leq E^{q-1} (T + at)^{\alpha (m-1) +1}+(T+at)  \; ,
\]
where $k = \min \{1-T^{- \alpha (m-q)} , 1-T^{- \alpha (1-m)}\}$. As it is simple to verify, the previous inequality holds for sufficiently small $a$. 

The proof of the second part is similar. In this case we choose $H(x,t)=(T-at)^{\alpha} E(x)$ for $q<1$, again with $T<1$. The only difference between $H$ and $G$ is to replace $\alpha$ by $-\alpha$. Consequently, we must show that:
\[
-a\alpha(T-at)^{\alpha-1}E\geq (T-at)^{\alpha m}(E^{m-1}E')'+(T-at)^{\alpha p}E^p -(T-at)^{\alpha q}E^q  \; .
\]
As a result, we find that: 

\[
a\alpha(T-at)^{\alpha-1}E\leq E^q ((T-at)^{\alpha q}-(T-at)^{\alpha m})+
E^p ((T-at)^{\alpha m}-(T-at)^{\alpha p})  \; .
\]

What remains is entirely similar to the previous case (less involved, actually) since we can work only with the first summand of the previous inequality (in the previous paragraph this was also feasible, but would restrict the proof to the case $m>1$). For the present situation we choose $\alpha = (p-q)(p-m)^{-1}(1-q)^{-1}$ if $p>m$ while $\alpha = 1/(m-1)$ if $p=m$, with $T=(\max_{x\in S}\{u_0(x)/E(x)\})^{1/\alpha}$ where $S=\{x:u_0(x)>0\}$. In addition, for $q=1$ we take $H(x,t)=(T+ a t)^{- \alpha} E(x)$, $T>1$, $\alpha = (m-1)^{-1}$ and the rest of the proof is similar to the one in the first part for the case $q=1$.$\:\;\:\; \Box$

\mbox{}

\mbox{}

\begin{flushleft}
{\bf 4. A criterion for the general case }
\end{flushleft}

The following result involving blow-up and extinction is more general than Theorem 1 in the sense that no parameter restrictions are imposed in the framework of problem (\ref{p1}).

\mbox{}

\noindent{\bf Theorem 2.}
{\em 
Let $u(x,t)$ be a solution of problem (\ref{p1}) with $p>q>0$, $p>1$, $m>0$. We then have: 
\begin{enumerate}
\item There exist initial conditions $u_0(x)\in C\cap L_1$ such that $u(x,t)$ presents blow-up.
\item If $u_0(x)\in C\cap L_1$ and $||u_0(x)||_\infty \leq 1$, then $u(x,t)$ presents either extinction or vanishing in an infinite time.
\end{enumerate}
}

\mbox{}

\noindent{\bf Proof.}

Item 2 can be shown by directly checking that if $||u_0(x)||_\infty \leq 1$, then the solution of the problem containing only the diffusion term (porous media equation) is a supersolution of problem (\ref{ec1}). As it is well-known (for instance, see \cite{Vaz}) the solutions of the porous media problem satisfy that: $\displaystyle{\lim_{t \to \infty}||u(x,t)||_\infty=0}$.

In order to prove item 1, let us show that:
\begin{equation}
\label{sub}
    \underline{u}(x,t)=(T-t)^{-\alpha} A(1-(\xi/a)^2)_+^b  \; ,
\end{equation}
where $(x)_+  =\{x \mbox{ if }  x>0; \mbox{ and } 0 \mbox{ if } x\leq 0\}$, $\xi=x(T-t)^{-\beta}$, $\alpha=1/(p-1)$ and $\beta=\alpha(p-m)/2$ is a subsolution of (\ref{ec1}). Showing (\ref{sub}) implies that every solution for which the initial condition $u_0(x)$ verifies 
$u_0(x) \geq \underline{u}_0(x)=T^{-\alpha} A(1-(xT^{-\beta})^2/a^2)_+^b$ leads to blow-up (in finite time). After substitution of (\ref{sub}) in (\ref{ec1}) we obtain:
\[
\begin{array}{ll}\alpha(1-(\xi/a)^2)_+^b-2b\beta(\xi/a)^2(1-(\xi/a)^2)_+^{b-1}\leq
4\frac{A^{m-1}}{a^2}b(mb-1)(\xi/a)^2(1-(\xi/a)^2)_+^{mb-2}\\
-2\frac{A^{m-1}}{a^2}b(1-(\xi/a)^2)_+^{mb-1}+A^{p-1}(1-(\xi/a)^2)_+^{bp}
-(T-t)^{\alpha(p-q)}A^{q-1}(1-(\xi/a)^2)_+^{bq} \; . \end{array}  
\]
Now let $I_i \;,\; i=1,\ldots,6$ be the 6 terms in the previous inequality, in their respective orders. Then such inequality can be rewritten in the form: 
\[
I_1-I_2\leq I_3-I_4+I_5-I_6 \Rightarrow I_2-I_1+I_3-I_4+I_5-I_6\geq0  \; .
\]
In order to prove the previous inequality we shall show that:
\begin{eqnarray}
  I_3/4+I_5/4-I_1 &\geq&0  \; ,  \label{des11} \\
  I_3/4+I_5/4-I_4 &\geq& 0 \; , \label{des12} \\
  I_3/4+I_5/4-I_6 &\geq& 0 \; , \label{des13} \\
  I_2+I_3/4 &\geq& 0       \; . \label{des14}
\end{eqnarray}
Note that with the choice $b>1/m$ all the $I_i$ are positive, excepting $I_2$ which depends on the sign of $\beta$. We then have: 
\begin{itemize}
  \item Regarding (\ref{des11}), let us first see the values of $\xi$ for which $I_3/4\geq I_1$. We want to prove that 
$$Bb(mb-1)(\xi/a)^2(1-(\xi/a)^2)_+^{mb-2}\geq\alpha(1-(\xi/a)^2)_+^b  \; ,$$ 
which implies:
$$Bb(mb-1)(\xi/a)^2\geq\alpha(1-(\xi/a)^2)_+^{2-b(m-1)} \; ,$$ where $B=\frac{A^{m-1}}{a^2}$. 
Assuming $b(m-1)<2$, and taking into account that $(1-(\xi/a)^2)_+\leq1$, the previous inequality holds if $(\xi/a)^2\geq \alpha /Bb(mb-1)=C_1$. Let us now verify when $I_5/4\geq I_1$, namely 
      $$A^{p-1}(1-(\xi/a)^2)_+^{bp}\geq 4\alpha(1-(\xi/a)^2)_+^b\Rightarrow (\xi/a)^2\leq 1-(4\alpha)^{\frac{1}{p-1}}A^{-\frac1b}=C_2 \; .$$
Choosing $a^2=A^{m-1-(p-1)/2}$, and thus $B=A^{(p-1)/2}$, we can now select a sufficiently large $A$ such that $C_2>C_1$, and this proves inequality (\ref{des11}).

 \item For (\ref{des12}), following a similar analysis to the one in \textit{(i)} we shall first work with $I_3/4\geq I_4$. Symplifying we obtain:

 $$(mb-1)(\xi/a)^2\geq2(1-(\xi/a)^2)_+\Rightarrow (\xi/a)^2\geq 2/(mb+1)=C_3 \; .$$
Note that $2/(mb+1)<1$. Let us now see when $I_5/4\geq I_4$:

 $$A^{p-1}(1-(\xi/a)^2)_+^{b(p-m)+1}\geq 2bB \Rightarrow (\xi/a)^2\leq1-(2bA^{-\frac{p-1}2})^{\frac{1}{b(p-m)+1}}=C_4 \; .$$

Therefore, if we choose a sufficiently large $A$, then $C_4>C_3$ and (\ref{des12}) is proven. Here we have taken into account the value of $B$ as a function of $A$ in the previous item, and we have assumed that $b(p-m)+1>0$. If $b(p-m)+1\leq0$ and taking a sufficienly large $A$ the inequality $I_5/4\geq I_4$ holds for all $(\xi/a)^2\geq0$.

\item For (\ref{des13}), we begin with $I_5/4\geq I_6$, from which we find
$$A^{p-q}(1-(\xi/a)^2)_+^{b(p-q)} \geq 4(T-t)^{\alpha(p-q)}\Rightarrow (\xi/a)^2\leq1-4^{\frac{1}{b(p-q)}}A^{-\frac1b}(T-t)^{\frac\alpha b} \; .$$
Consequently the last relationship will be fulfilled provided
$$(\xi/a)^2\leq 1-4^{\frac{1}{b(p-q)}}A^{-\frac1b}T^{\frac\alpha b}=C_5 \; .$$
Following a similar analysis to the one used with $I_3/4\geq I_6$, we have:
$$Bb(mb-1)(\xi/a)^2\geq(T-t)^{\alpha(p-q)}A^{q-1}(1-(\xi/a)^2)_+^{b(q-m)+2} \; .$$ 
Since $(1-(\xi/a)^2)_+<1$ and assuming $b(q-m)+2>0$, for the verification of the previous inequality it is sufficient that
 $$(\xi/a)^2\geq(T-t)^{\alpha(p-q)}A^{q-1}/Bb(mb-1) \; .$$
The previous condition is verified if
 $$(\xi/a)^2 \geq T^{\alpha(p-q)}A^{q-1}/Bb(mb-1)=C_6 \; .$$
Now taking $T=A^{-\frac{(p-1)^2}{2(p-q)}}$ we see that $C_6\sim A^{-(p-q)}$, so that taking a sufficiently large $A$ we have that $C_5>C_6$ and this shows (\ref{des13}).

 \item Regarding (\ref{des14}), if $b<1/(m-1)$ and $A$ is sufficiently large, this inequality always holds.
\end{itemize}
To conclude, let $\kappa = \max \{ m-q,m-1,0 \}$. Then note that if $1/m < b < 1/ \kappa$ (in particular, $1/m < b$ if $\kappa =0$) then the conditions imposed to $b$ in the proof shall be met.

The proof of Theorem 2 is now complete. $\:\;\:\; \Box$

\mbox{}

As a direct consequence of the results presented, it is worth noting that there exist initial conditions $u_0(x)$ leading to three different types of asymptotic behaviors: (i) blow-up solutions; (ii) vanishing solutions; and (iii) solutions for which $\displaystyle{\lim_{t\to \infty}||u(x,t)||_\infty=k\ne0}$.

\mbox{}

\mbox{}

\pagebreak
\begin{flushleft}
{\bf 5. Examples}
\end{flushleft}

\mbox{}

\noindent {\em 5.1 Example 1: Previously known solution family as a particular case}

\mbox{}

A particular case of the results displayed in Section 3 was considered by Samarskii and Galaktionov \cite{sgk}. In their analysis, the case $m=p>1$, $q=1$ is considered and explicit solutions of equation (\ref{ec1}) are found by means of the separation of variables technique, leading to solutions of the form $u_1(x,t)=\phi(t)\theta(x)$ where
\[
\phi(t)=e^{-t}\left(\frac{e^{-(m-1)t}}{m-1}+C_0 \right )^{-1/(m-1)} \; ,
\]
and 
\[
\theta(x)=\left\{
                \begin{array}{ll}
                  \left(\frac{2m}{(m^2-1)}\cos^2(\frac{\pi x}{L}) \right )^{1/(m-1)}, & |x|<L/2  \; , \\
                  0, & |x|\geq L/2\quad (L=\frac{2\pi\sqrt{m}}{m-1})  \; .
                \end{array}
              \right.  
\]
Depending on the sign of $C_0$ it can be seen that different behaviors can be found for $u_1$: it is stationary if $C_0=0$; it vanishes when $t\to\infty$ if $C_0>0$; and it blows-up when $C_0<0$. Now this result can be regarded as a particular case of Theorem 1. 

\mbox{}

\noindent {\em 5.2 Example 2: Further application of the separation of variables method}

\mbox{}

In a way similar to the one in the previous example, an additional family of explicit solutions can be found by means of the separation of variables method, namely $u_2(x,t)=\phi(t)\theta(x)$ for $q=m$, $p=1$. After substitution of the expression for $u_2$ in (\ref{ec1}) we arrive at:
\[
\frac{\phi'(t)-\phi(t)}{\phi^m(t)}=\frac{(\theta^{m-1}(x)\theta'(x))'-\theta^m(x)}{\theta(x)}=\lambda  \; .
\]
For convenience let us set $\lambda =-1$. Then it is not difficult to evaluate $\phi(t)$ in order to find: 
\[
\phi(t)=e^t(e^{-(1-m)t}-C)^{1/(1-m)}  \; .
\]
In order to evaluate $\theta(x)$ it is only necessary to recall that it verifies the stationary equation of (\ref{ec1}) for $q=m$, $p=1$. Such equation has an explicit solution which was constructed in Section 2, leading to:
\[
\theta(x)=\left(\frac{m+1}{2m}(1-\tanh^2(k_1 x/2))\right )^{1/(1-m)}  \; ,
\]
where $k_1=(1-m)/\sqrt{m}$. It can be seen that if $C<0$ then $\phi(t)\to\infty$ for $t\to\infty$; in second term, if $C>0$ then we have $\phi(t)\to 0$ when $t\to T=(m-1)^{-1} \log(C)$; and finally it is stationary if $C=0$.

\mbox{}

\noindent {\em 5.3 Example 3: A numerical illustration}

\mbox{}

For the sake of illustration, we shall now present a purely numerical example in which the solution $u$ either blows-up (Figure 1) or vanishes (Figure 2) depending on the initial condition to be greater (in the case of blow-up) or smaller (in the case of vanishing) than the stationary solution. To see this, we shall take $m=p$ because the stationary solution is explicit in this case, and for a better display we set the values $m=p=2$ and $q=0.9$. As initial conditions we have chosen twice the stationary solution (blow-up case) and half the stationary solution (vanishing case), where such stationary solution corresponds to instance (b) constructed in Section 2. The behaviors predicted by Theorem 1 are clearly visible in this case. 

\begin{center}
\fbox{\bf Figure 1 and Figure 2 here}
\end{center}

\mbox{}

\mbox{}

\begin{flushleft}
{\bf 6. Final remarks}
\end{flushleft}

In the Introduction we have seen that equation (\ref{ec0}) describes a wide variety of non-linear problems in different physical domains such as reaction-diffusion or growth dynamics. As it has been mentioned, the previous literature shows that this has been the case for both complementary possibilities $q>p$ and $p>q$. 

In this work, new results have been developed which are analogous in this context to the well-known existence of separatrices in the case of dynamical systems of the ODE type, a classical subject having clear consequences involving the properties and behavior of the solutions, including fundamental aspects such as stability, boundedness, etc. In our PDE case (\ref{ec0}), the parallelism is clear as far as the separatrix-like solution disjoints either blow-up (thus representing an unstable dynamical regime) or growth solutions (relevant in many applications, as seen in Section 1) from those solutions showing extinction (recall the classification in Definition 2). Moreover, the developments found sometimes allow a constructive approach, which actually leads to the generalization of previously known solutions and the construction of additional new ones, as illustrated in the Examples section. Consequently, we believe that the previous contributions have a clear interest in the analysis of problem (\ref{ec0}) and its physical interpretation.

\pagebreak

\pagebreak

\noindent {\bf FIGURE CAPTIONS}

\mbox{}

\mbox{}

\noindent {\bf Figure 1.} Numerical solution of equation (\ref{ec1}) for parameters $m=p=2$, $q=0.9$ and initial condition twice the stationary solution leading to blow-up (see Example 3 in Subection 5.3 for details).

\mbox{}

\noindent {\bf Figure 2.} Numerical solution of equation (\ref{ec1}) for parameters $m=p=2$, $q=0.9$ and initial condition half the stationary solution leading to extinction (see Example 3 in Subection 5.3 for details).


\begin{thebibliography}{99}
\bibitem{fis} R.A. Fisher, The wave of advance of advantageous genes, {\em Ann. Eugen.} 7 	(1937) 355-369.
\bibitem{kpp} A. Kolmogoroff, I. Petrovsky and N. Piscounoff, \'Etude de l'\'equation de 	la 	diffusion avec croissance de la quantit\'e de mati\`ere et son application \`a un 	probleme biologique, {\em Bull. Univ. Moskov}, Ser. Internat. Sec. Math. 1 (1937) 1-25.
\bibitem{arw} D.G. Aronson and H.F. Weinberger, Nonlinear diffusion in population genetics, 	combustion,and nerve pulse propagation, in: {\em Partial Differential Equations and 	Related Topics} (Lecture Notes in Mathematics Vol. 446, Springer-Verlag, Berlin, 1975), 	pp. 5-49.
\bibitem{GK} B.H. Gilding and R. Kersner, {\em Travelling waves in nonlinear
	diffusion-convection-reaction\/} (Memorandum No. 1585), Faculty of Mathematical 
	Sciences, University of Twente (2001).
\bibitem{PV} A. de Pablo and J.L. V\'{a}zquez, Travelling waves and finite propagation in a
	reaction-diffusion equation {\em J. Differential Equations\/} 93 (1991) 19-61.
\bibitem{SH} A. S\'{a}nchez-Vald\'{e}s, B. Hern\'{a}ndez-Bermejo, New travelling wave 
	solutions for the Fisher-KPP equation with general exponents, {\em Appl. Math. Lett.} 
	18 (2005) 1281-1285.
\bibitem{slh} M.J. Simpson, K.A. Landman, B.D. Hughes, D.F. Newgreen, Looking inside 
	an invasion wave of cells using continuum models: Proliferation is the 
	key, {\em J. Theoretical Biol.} 243 (2006) 343-360.
\bibitem{man} M.B.A. Mansour, Traveling wave patterns in nonlinear reaction-diffusion 	equations, {\em J. Math. Chem.} 48 (2010) 558-565.
\bibitem{egs} R. Engui\c{c}a, A. Gavioli, L. S\'{a}nchez, A class of singular first-order 	differential equations with applications in reaction-diffusion, {\em Discrete and 	Continuous Dynamical Systems} 33 (2013) 173-191. 
\bibitem{VO} E.O. Voit, {\em Computational Analysis of Biochemical Systems: A Practical Guide 
	for Biochemists and Molecular Biologists,\/} Cambridge University Press, Cambridge UK
	(2000).
\bibitem{M} J.D. Murray, {\em Mathematical Biology,\/} Springer-Verlag, New York (1993).
\bibitem{banks} R.B. Banks, {\em Growth and Diffusion Phenomena,\/} Springer-Verlag, Berlin (1994).
\bibitem{vos} E.O. Voit, M.A. Savageau, Analytical Solutions to a Generalized Growth Equation, {\em J. Math. Anal. 
Appl.} 103 (1984) 380-386.
\bibitem{marus1} M. Maru\u{s}i\'{c}, Z. Bajzer, Generalized Two-Parameter Equation of Growth, {\em J. Math. Anal. Appl.} 179 (1993) 446-462.
\bibitem{marus2} Z. Bajzer, M. Maru\u{s}i\'{c}, S. Vuk-Pavlovi\'{c}, Conceptual Frameworks for Mathematical Modeling of Tumor Growth Dynamics, {\em Math. Comput. Modeling} 23 (1996) 31-46.
\bibitem{ker} R. Kersner, Degenerate parabolic equations with general nonlinearities, {\em Nonlinear Anal.} 4 (1980) 1043-1062.
\bibitem{sgk} A.A. Samarskii, V.A. Galaktionov, S.P. Kurdyumov, A.P. Mikhailov, {\em Blow-up in 	Problems for Quasilinear Parabolic Equations,\/} (Nauka, Moscow, 1987 (in Russian); 	English translation: Walter de Gruyter, Berlin, 1995).
\bibitem{DL} K. Deng, H.A. Levine, The role of critical exponents in blow-up theorems: the 	sequel, {\em J. Math. Anal. 	Appl.} 243 (2000) 85-126.
\bibitem{Ka2} A.S. Kalashnikov, Some problems of the qualitative theory of second-order 	nonlinear degenerate parabolic equations (in Russian) {\em Uspekhi Mat. Nauk.} 42 (1987), 
	no. 2(254), 135-176. 
\bibitem{Ki} J.R. King, Extremely High Concentration Dopant Diffusion in Silicon, {\em IMA J. 	Appl. Math.} 40 (1988) 163-181.
\bibitem{GV3} V.A. Galaktionov, J.L. V\'{a}zquez, The problem of blow-up in nonlinear 	parabolic equations. Current developments in partial differential equations, {\em Discrete Contin. Dyn. Syst.} 8 (2002) 399-433.
\bibitem{VH} A.I. Vol'pert, S.I. Hudjaev, {\em Analysis in classes of discontinuous functions 	and equations of mathematical physics}, Mechanics: Analysis, 8. Martinus Nijhoff 	Publishers, Dordrecht (1985).
\bibitem{K2} R. Kersner, Nonlinear heat conduction with absorption: space localization and 	extinction in finite time, {\em SIAM J. Appl. Math.} 43 (1983) 1274-1285.
\bibitem{HVe} M.A. Herrero, J.J.L. Vel\'{a}zquez, Approaching an extinction point in 
	one-dimensional semilinear heat equations with strong absorption, {\em J. Math. Anal. 	Appl.} 170 (1992) 353-381.
\bibitem{FV} R. Ferreira, J.L. V\'{a}zquez, Extinction behaviour for fast diffusion equations with absorption, {\em Nonlinear Anal. Ser. A: Theory Methods} 43 (2001) 943-985.
\bibitem{eth} P. \'{E}rdi, J. T\'{o}th, {\em Mathematical Models of Chemical Reactions,\/} 	Manchester University Press, Manchester UK (1989).
\bibitem{bor} M. Borelli, M. Ughi, The fast diffusion equation with strong absorption: the 	instantaneous shrinking phenomenon, {\em Rend. Istit. Mat. Univ. Trieste} 26 (1994) 
	109-140.
\bibitem{Ka1} A.S. Kalashnikov, The effect of absorption on heat propagation in a medium 	in 	which the thermal conductivity depends on temperature, {\em Zh. Vychisl. Mat. Mat. Fiz.} 	16 (1976) 689-696.
\bibitem{abr} M. Abramowitz and I.A. Stegun, {\em Handbook of Mathematical Functions,\/} Dover, New York (1972).
\bibitem{arf} G. Arfken, {\em Mathematical Methods for Physicists,\/} Third Edition, Academic Press, Orlando FL  (1985).
\bibitem{Vaz} J.L. V\'azquez, {\em The Porous Medium Equation, Mathematical Theory,\/} Oxford 
	University Press, Oxford UK (2007).
\end{thebibliography}
\end{document}